\documentclass[12pt]{article}
\usepackage{epsfig}
\usepackage{epsf}
\usepackage{here}
\usepackage{color}
\usepackage{pifont}

\newcommand{\DMdet}{\mbox{\boldmath ${\cal D}$}}
\newcommand{\1}[1]{1_{\{#1\}}}
\newcommand{\idp}[1]{\perp \hspace{-0.25cm} \perp_{#1}}
\newcommand{\infl}[1]{\rightarrow\rightarrow_{#1}}
\newcommand{\dinf}[1]{\longrightarrow_{#1}}
\newcommand{\sli}[1]{\mbox{$\rightarrow \rightarrow\hspace{-0.60cm} / \hspace{0.40cm}_{#1}$}}
\newcommand{\wli}[1]{\mbox{$\longrightarrow \hspace{-0.60cm} / \hspace{0.40cm}_{#1}~$}}

\newcommand{\FF}{{{\cal F}}}
\newcommand{\HH}{{{\cal H}}}
\newcommand{\Ob}{{{\cal O}}}
\newcommand{\X}{{{\cal X}}}

\newcommand{\LL}{{\cal L}}

\newcommand{\A}{{\cal A}}
\newcommand{\PS}{{\cal S}}

\newcommand{\ind}{\perp \hspace{-0.25cm} \perp}

\newcommand{\EPZ}{{\rm E}_{P_0}}

\newcommand{\bX}{\mbox{\boldmath $X$}}
\newcommand{\bA}{\mbox{\boldmath $A$}}
\newcommand{\bS}{\mbox{\boldmath ${\cal S}$}}

\newtheorem{Theorem}{Theorem}

\newtheorem{Lemma}{Lemma}
\newtheorem{Proposition}{Proposition}
\newtheorem{Definition}{Definition}
\begin{document}

\title{A general dynamical statistical model with possible causal interpretation}
\author{ Daniel Commenges$^{1,2}$ and Anne G\'egout-Petit $^{2,3}$ }
\maketitle
\noindent {\em 1  INSERM, U 875, Bordeaux,  F33076, France }\\
{\em 2 Universit\'e Victor Segalen Bordeaux 2, Bordeaux, F33076, France}\\
{\em 3 IMB, UMR 5251, Talence, F33405, France}\\
%Tel: (33) 5 57 57 11 82; Fax (33) 5 56 24 00 81; email daniel.commenges@isped.u-bordeaux2.fr

\maketitle {\bf Summary}. We develop a general dynamical model as a framework for possible causal
interpretation. We first state a criterion of local independence in terms of measurability of processes
involved in the Doob-Meyer decomposition of stochastic processes, as in Aalen (1987); then we
define direct and indirect influence. We propose a definition of causal influence using the concepts of
``physical system''. This framework makes it possible to
link descriptive and explicative statistical models, and encompasses quantitative processes and events.
One of the features of this paper is the clear distinction between the model for the system and the model for the
observation. We give a dynamical
representation of a conventional joint model for HIV load and CD4 counts. We show its inadequacy to capture causal influences while on the contrary known mechanisms of HIV infection can be expressed directly through a system of differential
equations.\vspace{3mm}

{\em Keywords}: Causality; causal influence; differential equations; directed graphs; dynamical models; HIV; randomisation; stochastic processes.

\section{Introduction}
The issue of causality has been studied by many philosophers since Aristotle and is of central importance in all
 branches of science (see Bunge, 1979 and Salmon, 1984). A central question for scientists who
use statistics and for statisticians is whether statistical models may help in deciphering causal links. After
recognising that correlation is not causation, scientists have tended to use statistical methods as one element
among others to help establish causal links. Epidemiologists are particularly cautious, and with good
reason, in concluding to causal influences. There has been however a growing interest in developing statistical
models able to represent causal influences. From the beginning, graphs have played an important role in
representing the set of causal influences. The pioneering work of Wright (1921, 1934) have inspired the more
recent developments of structural equation models (Joreskog, 1978) and graphical models (Dawid, 1979; Lauritzen
and Wermuth, 1989; Cox and Wermuth, 1996). An approach using the modelling of ``potential outcome'', often called
the counterfactual approach, has been proposed in the context of clinical trials by Rubin (1974) and further
studied by Holland (1986) among others. The counterfactual approach has been extended to the study of
longitudinal incomplete data in several papers, the results of which have been gathered together by van der
Laan and Robins (2002). Spirtes, Glymour and Sheines (2000) and Pearl (2000) develop the issue of investigating
causality with graphical models.

The counterfactual approach however has been criticised (Dawid, 2000; Geneletti, 2007) and the modelling of potential outcomes
raises difficulties when treating truly dynamical problems. In fact another school tackles causality by
directly using dynamical models. This approach started in the econometrics literature with Granger (1969) and Schweder
(1970) and was more recently developed by several Scandinavian statisticians using the formalism of stochastic
processes, and in particular of counting processes (for a review see Eerola, 1994; Aalen and Frigessi, 2007). Of
particular interest is the paper by Aalen (1987) which outlines a general approach for defining influences for
stochastic processes through the Doob-Meyer decomposition. The most recent developments of the dynamical approach
are the method of ``dynamic path analysis'' of Fosen et al. (2006)  and the study of the possibly cyclic
directed graphs associated with this definition of influence by Didelez (2007). Defining influence in the stochastic process framework does not ensure that we make relevant causal
inference but we believe that it provides a better formalism for tackling this issue than
approaches which deal only with random variables.

The aim of this paper is to develop the dynamical approach in a general framework, focusing in particular on causal interpretation, using the concept of system, which was advocated long ago by von Bertallanffy (1968); we attempt to go from the concept of ``influence'', which is mathematically defined, to the concept of ``causal influence'', which has a physical meaning. We make a clear distinction between the  model for the system and the model for the observations, a classical distinction in automatics (Jazwinsky, 1970) but not in biostatistics. Moreover we link classical epidemiological models and mechanistic models; the latter are not generally taken into consideration in the literature of causal models although (or because) they make explicit use of scientific knowledge.

The paper is organised as follows.
In section 2 we develop a criterion of local independence in terms of measurability of processes involved in the Doob-Meyer representation; then we define direct and indirect influence. In section 3 we propose a definition of causal influence using the concepts of ``physical system'' and ``physical laws'' for which we propose a definition. Our framework makes it possible to link descriptive and explicative statistical models and encompasses the analysis of events and of quantitative processes. In section 4 we develop the distinction between the model for the system and the model for the observation.
In section 5 descriptive and explicative joint models of HIV load and CD4 counts are considered.

\section{Local independence, direct and indirect influence}
\label{infl}
\subsection{Notations}
Consider a filtered space $(\Omega, \FF, (\FF_t), P)$ and a multivariate stochastic process $\bX=(\bX_t)_{t\ge 0}$; $\bX_t$ takes  values in $\Re^m$, and the whole process $\bX$ takes values in $D(\Re^m)$, the Skorohod space of all cadlag functions: $\Re_+ \rightarrow \Re^m$.  We suppose that all the filtrations satisfy the usual conditions. We have $\bX=(X_j, j= 1, \ldots,m)$ where $X_j=(X_{jt})_{t\ge 0}$. We shall note $X_j \in \bX$. We denote by $\X_t$ the history of $\bX$ up to time $t$, that is  $\X_t$ is the $\sigma$-field $\sigma(\bX_u, 0\le u \le t)$, and by $(\X_t)=(\X_t)_{t\ge 0}$ the families of these histories, that is the filtration generated by $\bX$. Similarly we shall denote by $\X_{jt}$ and $(\X_{jt})$ the histories and filtration associated to $X_j$. If $C$ is a subset of $(1,\ldots,m)$ we shall call $X_C$ the multivariate process $(X_j, j\in C)$.

\subsection{Local independence, direct and indirect influence}
Let $\FF_t=\HH \vee \X_t$; $\HH$ may contain information known at $t=0$, in addition to the initial value of $\bX$. We shall consider the class of special semi-martingales, that is the class of processes which admit a unique Doob-Meyer decomposition in the $(\FF_t)$ filtration, under probability $P$:
\begin{equation} \bX_t=\Lambda_t + M_t, t\ge 0, \end{equation}
where $M_t$ is a martingale and $\Lambda_t$ is a predictable process with bounded variation. We shall denote the Doob-Meyer decomposition of $X_j$: $X_{jt}=\Lambda_{jt} + M_{jt}$.
We shall consider the non-degenerate case in which all the components of $M$ are different from zero; the
deterministic case will be studied in section \ref{detcase}. We shall assume two conditions bearing on the
bracket process of the martingale $M$:

{\bf A1} $M_{j}$ and $M_{k}$ are orthogonal martingales, for all $j\ne k$;

{\bf A2} $X_j$ is either a counting process or is continuous with a deterministic bracket process, for all $j$.

We call $\DMdet$ the class of all special semi-martingales satisfying {\bf A1} and {\bf A2}. The class of special semi-martingales is stable by change of absolutely continuous probability (Jacod and Shiryaev, 1987, page 43) and this is also true for the the class $\DMdet$.

\begin{Definition}[Weak conditional local independence (WCLI)]

$X_k$ is weakly locally independent of $X_j$ in $\bX$ on $[0,\tau]$ if and only if $ \Lambda _k$ is
$(\FF_{-jt})$-predictable on $[0,\tau]$, where $\FF_{-jt}=\HH \vee \X_{-jt}$ and $\X_{-jt}=\vee _{l\ne
j}\X_{-lt}$. Equivalently we can say in that case that $X_k$ has the same Doob-Meyer decomposition in
$(\FF_{t})$ and in $(\FF_{-jt})$. We will note in that case $X_j \wli{\bX} X_k$.
\end{Definition}

{\bf Remark 1}. Assumption {\bf A2} is necessary for the measurability-based definition of WCLI to be clearly
interpreted. If we did not impose {\bf A2} we could find counter-examples in which a WCLI holds while intuitively
independence does not hold. Such a counter-example is the process $\bX=(X_1, X_2)$ which is the solution of the
differential equation: $dX_{1t}=a~dt+b~dW_{1t}$; $dX_{2t}=X_{1t}~dt+e^{X_{1t}}~dW_{2t}$, Where $W_1$ and $W_2$ are Brownian motions. We would not like to
say that $X_2$ is WCLI of $X_1$. However, because $X_1$ appears in the bracket process of $X_2$, $\X_{1t}$ is
included in $\X_{2t}$ so that $ \Lambda _2$ is $\X_{2t}$-predictable and thus we would conclude that $X_2$ is
WCLI of $X_1$.

{\bf Remark 2}.
It is tempting to define WCLI directly in terms of the conditional independence:
\begin{equation} \X_{kt} \idp{\X_{Ct-},\X_{kt-}} \X_{jt-} , 0\le s<t\le \tau.\end{equation}
Here $\bX=(X_j,X_k,X_C)$. 
However, this condition is void in general when we consider processes in continuous time. Because conditional independence is defined via conditional probability and in general, events of $\X_{kt}$ will have conditional probabilities equal to one or zero given  $\X_{kt-}$, the condition will always hold. It is possible that WCLI can be defined in terms of conditional independence of $\sigma$-fields but this is an open problem.

\begin{Definition}[Direct influence]
We shall say that if $X_k$ is not WCLI of $X_j$ in $\bX$, $X_j$ directly influences $X_k$ in $\bX$ and we will note $X_j \dinf{\bX} X_k$.
\end{Definition}

\begin{Definition}[WCLI and Direct influence for set of components]
Let $A, B$ subsets of  $(1,\ldots,m)$.
We shall say that  $X_A \dinf{\bX} X_B$ if there is $j \in A$ and $k \in B$ such that $X_j \dinf{\bX} X_k$.
\end{Definition}

What we call here ``direct influence'' is the time-continuous analogue of Granger strong causality (Granger,
1969). 
We may consider another, stronger, condition of local independence.

\begin{Definition}[Strong conditional local independence (SCLI)]
$X_k$ is SCLI of $X_j$ in $\bX$ if and only if $X_j\wli{\bX} X_k$ and
there is no $X_D \in \bX$ such that $X_j\dinf{\bX} X_D$ and $X_D \dinf{\bX} X_k$
and we will note $X_j \sli{\bX} X_k$
\end{Definition}

\begin{Definition}[Influence]
We shall say that if $X_k$ is not SCLI of $X_j$, $X_j$  influences (at least indirectly) $X_k$  in $\bX$  and we will note $X_j \infl{\bX} X_k$.
\end{Definition}

An interesting case is when weak independence holds but strong independence does not hold; equivalently
$X_j$ influences $X_k$ but $X_j$ does not directly influence $X_k$: we shall say that  $X_j$ indirectly
influences $X_k$.

\begin{Definition}[Indirect influence]
If $X_j \infl{\bX} X_k$ and $X_j \wli{\bX} X_k$ then there is $X_C \in \bX$ such that $X_j \dinf{\bX} X_C \dinf{\bX} X_k$ and we shall say that $X_j$ indirectly influences $X_k$ through $X_C$ in $\bX$.
\end{Definition}

{\bf Remark.} Since the Doob-Meyer decomposition depends on $P$ so do all the independencies and influences; realising this fact is crucial for the definition of causal influence in section \ref{causality-sect}.

\subsection{Differential equation: towards causal interpretation}
\label{ODE-sect}
Writing the process of interest in the form of a stochastic differential equation (SDE) is a way of making the causal mechanisms at work more explicit.
If $\Lambda_t$ is differentiable, the Doob-Meyer decomposition can be written:
\begin{equation} d\bX_t= \lambda_t dt+ dM_t, \end{equation}
with $\Lambda_t=\int_0^t \lambda_u du$.
Differential equation models are commonly used in physics, biology and in finance (Oksendal, 2000) to model the evolution of $\bX_t$ as a function of the past plus a random term brought by the martingale. The two main cases, which have been considered in different streams of research, are the case where the trajectories of $\bX$ are continuous and the case where $\bX$ is a counting process. In the case of continuous trajectories of $\bX$ it is common to take for $M$ a Brownian martingale (in which case $dM_t=f(t)dW_t$, with $W=(W_t)$ a Brownian motion), and the models considered are It\^o  processes.  In the case where $\bX$ is a counting process we write:
\begin{equation} d\bX_t= \lambda_t dt+ dM_t, \end{equation}
where $M$ is a discontinuous martingale with predictable variation process equal to $\Lambda$, and $\lambda$ is called the intensity of the process. We may consider mixing the two cases, considering that $\bX=(X_1, X_2)$, where $X_1$ is an It\^o process and $X_2$ a counting process, each of these processes being possibly multivariate.
The processes defined by these differential equations are not Markov in general. The Markov assumption is an interesting particular case and it is discussed in section \ref{Markov}.

\subsection{The deterministic case}
\label{detcase} Ordinary differential equation (ODE) models seem to arise as particular cases in which $M=0$.
So one way to apply our definition of WCLI to deterministic models is to consider that these models are in fact
stochastic but the martingale has a bracket process which takes small values in regard to $\Lambda$. Particular
phenomena appear in purely deterministic models, in particular because  the concept of {\em filtration} no longer applies. In
 that case the unicity of the differential equation is lost.
For instance consider the process $\bX=(X_1,X_2)$; consider the case where the process $\bX$ is deterministic
and the trajectories are solutions of the ODE system:  $dX_{1t}=a~dt$; $dX_{2t}=X_{1t}dt$ with initial
conditions $ X_{10}=X_{20}=0$. The trajectories are also solutions of the ODE system: $dX_{1t}=a~dt$;
$dX_{2t}=at~dt$  with initial conditions $X_{10}=X_{20}=0$. One would be tempted to say that $X_1$ influences
$X_2$ when looking at the first ODE system and but not when looking at the second one. The second ODE system
however is not time-homogeneous. Unicity can thus be restored if we impose the restriction of time-homogeneity
(see in section \ref{Markov} a discussion of the physical meaning of time-homogeneity). Taking advantage of the
unicity of the time-homogeneous differential equation representation, we will consider it as the canonical
representation, if it exists. We can then use the definition of WCLI for stochastic differential equations to
define WCLI for the deterministic case: construct a SDE system by adding to the ODE system orthogonal
martingales with deterministic brackets. The influence graph of the time-homogeneous ODE system is, by
definition, the same as that of the derived SDE. In the above example, if we add a standard Wiener martingale to
the canonical (time-homogeneous) representation we obtain the SDE: $dX_{1t}=a~dt+dW_{1t}$;
$dX_{2t}=X_{1t}dt+dW_{2t}$, in which it is clear that we have $X_1 \dinf{\bX} X_2$.

\subsection{Graph representation}
\begin{figure}[h]

\centering
\includegraphics[scale=0.6]{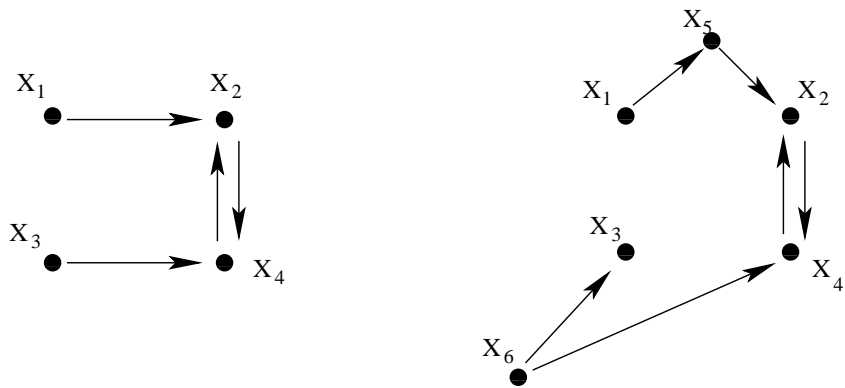}

\caption{Example of two graphs from the same physical system.\label{nested}}
\end{figure}

We may construct as in Didelez (2007) a directed graph representing influences between components of $\bX$. This directed graph has for vertices the components $X_j$ and there is a directed edge $(j,k)$ if and only if $X_j \dinf{\bX} X_k$. Note that there can be two directed edges between two vertices, for instance $(j,k)$ and $(k,j)$; this can be denoted by two arrows or by a double-sided arrow ($\longleftrightarrow$). A path is an ordered sequence of directed edges $\{(j_0, j_1), (j_1,j_2),\ldots,(j_{k-1},j_k)\}$. Indirect influence can be read directly off the graph: $X_j \infl{\bX} X_k$ if there is a path from $j$ to $k$. An example is shown in figure \ref{nested} which represents the hypothetical influence graphs for processes $\bX^1$ on the left and $\bX^2$ on the right; the graphs are not acyclic and in particular we have $X_2\dinf{\bX^1} X_4$ and $X_4\dinf{\bX^1} X_2$. We see also that $X_1$ indirectly influences $X_4$ but does not influence $X_3$, which we can note: $X_1\infl{\bX^1} X_4$ and $X_1\sli{\bX^1} X_3$. The graph for $\bX^2$ on the right may represent a richer system; we shall develop the issue of considering a family of nested systems in section \ref{causality-sect}.

 We say that $X_C$ blocks the paths from $X_l$ to $X_k$ if all the paths from $X_l$ to $X_k$ contain a node in $X_C$.
 For instance $X_4$ blocks the paths from $X_3$ to $X_2$ in $\bX^1$. In $\bX^2$ there is no path from $X_3$ to $X_2$, so $X_4$ still blocks the paths from $X_3$ to $X_2$, although in a trivial manner. If there is a path from $X_l$ to $X_k$ and $X_j$ blocks the paths from $X_l$ to $X_k$ there is necessarily a path from $X_l$ to $X_j$ and a path from $X_j$ to $X_k$, which can be expressed as:

\begin{Lemma}[Decomposable influence] \label{decomp-infl}
If $X_l \infl {\bX^m} X_k$ and $X_j$ blocks the paths from $X_l$ to $X_k$ then $X_l \infl {\bX^m} X_j$ and $X_j \infl {\bX^m} X_k$.
\end{Lemma}

\section{Causal influences}
\subsection{Systems, causal influence}
\label{causality-sect}
%We may postulate the existence of a detailed mechanistic graph generated by a multivariate process $X^*$ for a given system. We can consider a family of  graphs $(\G)$ where $\G$ is generated by a subset of components of $X^*$. Ideally $X^*$ is a homogeneous Markov model and the direct influences which appear in its graph can be interpreted as direct (or at least indirect) causal effect. This means that there is a stability of the dynamics. This means that adding components which are not influenced by $X^*$, but which may influence $X^*$ does not modify the dynamics of $X^*$. Pruning $X^*$ will often result in loss of the homNyogeneity and the Markov properties and may result in loss of some WCLI properties. The ``faithful'' property ensures that it will not make new WCLI properties appear.
In this section we outline a philosophical theory of causality; this theory is necessarily incomplete and questionable but we feel that a theory of this kind is necessary to link the mathematical definitions to the real world. The main concept used will be that of {\em system} and we will take two examples to illustrate this and other concepts: the first is the archetypal example of the solar system; the second is the system formed by the immune system and a population of HIV viruses. Thus our first task is to define a system, that we will also call
 a ``physical system'' $\PS$ in which we are interested. To define a system we admit that we can define a {\em level} at which relevant characteristics can be defined: we may distinguish a vector of {\em attributes} and a {\em state} vector. The attributes essentially define the system and do not vary in time; the state represents the characteristics, in general varying with time,  in which we are really interested and will be represented by a multivariate stochastic process. We may decide for instance that the level we are interested in is that of the sun and the planets and their trajectories. A possible system may be identified by the sun and the nine planets; the attributes of the system are the masses of these ten celestial bodies; the state at time $t$ is the vector of position and speed (in a reference system) of the ten celestial bodies.
What we have excluded in defining this level are the details of the planets such as their physical structure, presence of life, particular events like storms and so on (see Batterman, 2002).
 In the Immune-HIV system example, we may decide that the level we are interested in is that of populations of cells or of HIV viral particles in a particular subject; the attributes describe which types of cells or viral particles are considered and characteristics of these populations if they may differ from one subject to another; the state may be the numbers (or concentrations) of different types of cells or of viral particles. At this level we are not interested in the fate of a particular cell.

Now we suppose that we are particularly interested in one or several components of the state. We assume that there are laws which govern the evolution of the states and some of the laws tell us that the evolution of the component $j$ at time $t$ depends on the component $k$ just before $t$. Newton's laws, including the gravitation law, tell us how to compute the force of attraction between two massive bodies. However it is impossible to find a system with two massive bodies completely isolated from the rest of the universe;  moreover it is very difficult to avoid the circularity in the definitions. For instance we have used the concept of {\em system} in the previous sentence, a concept which is still not defined. It would be tempting to define {\em system} by first defining {\em causal influence}. However in order to define {\em causal influence} we need to apply natural laws to a system. There is the problem of defining {\em natural laws} or {\em physical laws}.

\begin{Definition} [System]
A system is the couple $\PS=(\bA,\bX)$ of attributes and state. The attribute $\bA$ is a possibly random element with value in $\Re ^d$ which, together with the state, is sufficient to identify the system.
The state $\bX$ is a stochastic process from $(\Omega, \FF)$ on $(D(\Re^m), \Sigma)$, where $\Omega$ is the universe and  $\FF$ contains all the events pertaining to the level of interest; $D(\Re^m)$ is a Skorohod space of all cadlag functions: $\Re_+ \rightarrow \Re^m$, and $\Sigma$ the Borel sigma-field derived from the Skorohod topology.
\end{Definition}

%There is an abuse of language in this definition in that the vector of attributes is really a set of qualities to which we make correspond a real value by the possibly random element $\bA$, and the same remark holds for the state. 
We consider that deterministic $\bA$ and $\bX$ is a particular case of the stochastic case, with the reservation made in section \ref{detcase} that WCLI is defined only for time-homogeneous ODE. Often the attribute will be considered as deterministic. In the solar system example, both attributes and states may be considered as deterministic, or we may consider it as random but work with a probability conditional to the observed value. The rationale for considering attributes as random is that they are the results of systems of another level: the existence of the planets is the result of the process of formation of the solar system. Note that even in that example, complex systems or long range predictions may raise the issue of {\em chaos}, thus introducing a stochastic feature (Murray and Dermott, 1999).

%In some applications it may be relevant to group components of attributes and states to form sub-systems or {\em entities}; for instance in the solar system example it would make sense to group the attributes and the components of the state of a given planet to form a sub-system or entity restricted to that planet; we would then restrict the attention to possible causal influences of one planet on another one. However we do not pursue this presentation in this paper because it does not make sense in most applications in epidemiology.

Given a system $\PS^{m}=(\bA^m, \bX^m)$, we call $\FF^m_t$ the
sigma-field generated by the attribute and the history of the state
at time $t$, $\FF^m_t= \sigma(\A^m) \vee \sigma(\bX^m_u, 0\le u \le
t)=\mathcal{A}^m \vee \X^m_t$. It is important to consider several
systems and in particular we may consider {\em nested systems}. A
system $\PS^{m'}$ is nested in $\PS^{m}$ if $\FF^{m'}_t\subset
\FF^{m}_t$ for all $t$ : $\PS^{m'}$ can be enlarged by addition of
attributes ($\mathcal{A}^{m'} \subset \mathcal{A}^m$) and/or
addition of $\bX^{m'}$ components  ($\X^{m'}_t \subset \X^m_t$). We
can consider a sequence of nested systems $\bS=\{\PS^{m}\}_{m>0}$
(we note $\PS^{m}\in \bS$ and  $\PS^{m}\subset \PS^{m'} $ if
$m < m'$). In this case, the family $\{\FF^{m}_t\}_{m>0}$ forms a
filtration (for each $t$). If we consider a period of observation
$[0,\tau]$ (included in the definition of the level) we note
$\FF^m=\FF^m_{\tau}$. Note that saying that $\PS^{m}\subset \PS^{m'} $ is more general than considering that all the components of $\PS^{m}$ belong to $\PS^{m'}$, although most result will refer to this case. 

From now on, we will speak about direct and indirect influences of
$X_j$ on $X_k$ in the system $\PS^m$ (and denote $X_j \dinf {\PS^m}
X_k$ or $X_j \infl {\PS^m} X_k$) these influences corresponds of the
definitions of influences in $\bX^m$ in section \ref{infl} with
${\cal H }=\A^m$.

We assume that there is a true probability law $P^*$ on $(\Omega,
\FF)$ and we denote its restriction to $\FF^m$ by $P^*_{\FF^m}$. We
would like to approach $P^*_{\FF^m}$ by applying {\em physical
laws}. We shall now endeavor to define {\em physical laws}. Let us
first define {\em mathematical laws}.

\begin{Definition} [Mathematical laws]
Mathematical laws at a certain level are a set of mathematical
procedures that can be applied to any system $\PS^{m}$ of this
level to build a probability $P^{\PS^{m}}$ on $\FF^m$.
\end{Definition}
Generally the probability $P^{\PS^{m}}$ will be different from
$P^*_{\FF^m}$. Suppose that we are particularly interested in a
system $\PS^{1}$, we may have to consider richer systems for making
correct predictions for the system of interest. We define {\em
physical laws} as yielding a probability that may be as close as we
wish from $P^*_{\FF^1}$, if we can apply them to a correct system.

\begin{Definition} [Physical laws] If for any system $\PS^{1}$ of a given level, there exists a sequence of nested systems $\bS=\{\PS^{m}\}_{m>0}$ including $\PS^{1}$ and  {\em mathematical laws} such that
 $P^{\PS^{m}}_{\FF^1}$ converges weakly toward $P^*_{\FF^1}$, these mathematical laws will be called {\em physical laws} at this level, and such a sequence $\bS$ will be called an approximating sequence for $\PS^{1}$.
\end{Definition}

The weak convergence means that $\int {g}~
dP^{\PS^{m}}_{\FF^1}\rightarrow \int {g}~dP^*_{\FF^1}$ for
any $\FF_1$-measurable continuous bounded function ${g}$  on
$\Omega$. We may also write $d_P(P^{\PS^{m}}_{\FF^1}, P^*_{\FF^1})
\rightarrow 0$, where $d_P(.,.)$ is the Prokorov metric for
probability measures based on the Skorohod topology. The advantage
of the Prokorov metric is that it metrizes weak convergence (Gibbs and Su, 2002) and it
encompasses the deterministic case (which makes sense in the solar
system example). In the deterministic case
$\bX$ takes the value $X^*$ with probability one under $P^*_{\FF^1}$
and the value $X^{\PS^{m}}$ under $P^{\PS^{m}}_{\FF^1}$ and we
have $d_P(P^{\PS^{m}}_{\FF^1},
P^*_{\FF^1})=d(X^{\PS^{m}},X^*)$.

We may postulate the existence of {\em physical laws}. This
postulate reflects the asymptotic {\em separability} of the
universe; that is, for making good predictions we do not need to
take into account the whole universe, but on the other hand,
application of the laws (even if we know the correct laws) never
leads to perfect prediction, partly because we have isolated a
system from the rest of the universe.

The systems may be more or less satisfactory according to the
distance to the true probability achieved. For instance we would not
call a set constituted of the Earth and Mars a satisfactory system;
if we applied Newton's laws to this system we would see that the
observed trajectories would be in large disagreement with the
predicted ones; we would thus search for a better set of bodies, for
instance the set (Sun, Earth, Mars).

We have to make an assumption of finiteness of the approximating
sequence to have a clear definition of causal influence. We
conjecture that this assumption could be avoided using a
quantitative approach of WCLI but this is beyond the scope of
this paper.

 {\bf A3}. There is a {\em perfect} system $\PS^{M}$ for
$\PS^{1}$ such that $\FF^1 \subset \FF^M$ and
$P^{\PS^{M}}_{\FF^1}=P^*_{\FF^1}$.

This means that the probability law computed with the physical law applied to system $\PS^{M}$ coincides with the true law on the events of interest $\FF^1$.

We assume that {\bf A1} and {\bf A2} hold for all the systems considered (see
discussion in section \ref{stability-properties}); assuming {\bf A3} we
can give the following definition.
\begin{Definition} [Causal influence] A component  $j$  has a causal influence on a component $k$ in $\PS^{1}$ if $X_j \infl {\PS^{M}} X_k$ under $P^*$, if
$\PS^{M}$ is a perfect system for $\PS^{1}$.
\end{Definition}

{\bf Remark}. The direct influences under the physical law are the
same in all the systems and in particular in the perfect system; a
direct influence under the physical law is thus always a causal
direct influence.

{\bf Example: Solar system}. If we consider a system (Earth, Moon)
the law of gravitation tells us that the earth (in our presentation,
the position of the Earth) has an influence on the trajectory of the
moon; by definition (if we accept that the law of gravitation is a {\em
physical law}), this is a causal influence. Even if this system is
not completely satisfactory, the notable fact is that in any richer
system, the Earth will have an influence on the Moon; this stability is
characteristic of causal influences.

{\bf Example: Immune system-HIV}. The mechanisms which derive from
the properties of HIV and CD4 lymphocytes are such that HIV can
infect CD4 lymphocytes and that infected lymphocytes can produce
viruses. The number of viruses produced depends in part on the
number of infected lymphocytes. Thus the component of the state
``number of viruses'' (in a given individual) has a causal influence
on the component ``number of infected lymphocytes''. We can deduce
the form of causal influences at the level of concentrations from
knowledge of the mechanisms which lead to the replication of the
virus and application of diffusion laws. The approach is similar to
Boltzman's theory of gases (see Strevens, 2005).

The problem of estimation of the true law $P^*$ will be dealt with
in section \ref{observations}.

\subsection{Stability of structures in sequences of systems}

\subsubsection{Stability of the class $\DMdet$}\label{stability-properties}
Since influences are defined in the class $\DMdet$ and there is a
need to consider sequences of systems, the stability of $\DMdet$ in
the sequence is crucial. We have the following Lemma:
\begin{Lemma}[Stability of $\DMdet$ ] \label{stability-class}
If $\bX^M \in \DMdet$, then $\bX^m\in \DMdet$ for all $\bX^m \subset
\bX^M$.
\end{Lemma}
    {\bf Proof.} The class of special semi-martingales is stable by change of filtration (Jacod and Shiryaev, 1987).
    The optional square-bracket process does not depend on the filtration. Thus it remains deterministic for continuous processes ({\bf A2})
    and the martingales remain orthogonal ({\bf A1}). For counting processes the orthogonality of the martingales holds if and only if the martingales cannot jump
    at the same time, which does not depend on the filtration; the martingales of a continuous and a counting process are always orthogonal.

\subsubsection{Instability of the homogeneous Markov property}
 \label{Markov}
It is interesting to examine the Markov properties of the models. In
the general case the derivatives of the predictable processes
involved in the Doob-Meyer decomposition depend on the whole past of
the process. For instance we can make these dependencies explicit by
writing $\lambda(t)= \lambda(t, X_{u}, 0\le u<t)$. In Markov models
these functions depend only on the present, or more precisely on
$X_{t-}$: $\lambda(t, X_{u}, 0\le u<t)= \lambda(t, X_{t-})$. The
model is (time)-homogeneous if these functions do not depend on
time: $\lambda(t, X_{t-})=\lambda(X_{t-})$. Typical physical models
are time-homogeneous Markov models. The Markov property means that
knowledge of the past cannot improve our knowledge of the future if
we know the present, and the homogeneous property means that the
laws of the universe we have used for constructing the model do not
change. So one can argue that if the model we consider is not a
homogeneous Markov model we have omitted important components in
the model.

 If a process is time-homogeneous Markov in $\PS^{M}$ under $P^*$ this does not in general hold for $\PS^{m} \subset \PS^M$; thus the homogeneous Markov property is not stable in $\DMdet$. This fact explains why it is often needed to consider non-homogeneous and even non-Markovian models in biology; indeed the systems considered are often oversimplified in view of the complexity of the real systems, leading to a loss of the time-homogeneous Markov property.
 
\subsubsection{Faithfulness and stability of influences}
 To go further, we assume that $P^*$ is ``faithful'', a property which is discussed for instance in Robins et al. (2003) for directed acyclic graphs, and that we define in our context as:

\begin{Definition}[Faithful probability]
 A probability $P$ is faithful for a sequence $\bS$ if for any $\PS^{m'},\PS^{m} \in \bS$ such that $\PS^{m'} \subset \PS^m$ and such that $X^m_j=X^{m'}_j=X_j$ and $X^m_k=X^{m'}_k=X_k$, we have $X_j \dinf {\PS^m} X_k$
  implies $X_j \dinf {\PS^{m'}} X_k$. Equivalently, $X_j \wli {\PS^{m'}} X_k$   implies $X_j \wli {\PS^m} X_k$.
\end{Definition}

Figure \ref{nested} illustrates a case which is compatible with $P$
being faithful: if $\PS^{1}$ (resp. $\PS^{2}$) has the left (resp.
right) influence graph, we see that the weak independence between $X_1$ and
$X_3$ is stable when the system is enriched from $\PS^1$ to $\PS^2$;
on the other hand the influence of $X_3$ on $X_4$ disappears ($X_6$
acts as a confounder process); finally the direct influence of $X_1$
on $X_2$ becomes an indirect influence through $X_5$. Faithfulness
does not hold in general; however one may argue that it does not
hold only in very specific cases. We show for instance in Appendix
A that

\begin{Proposition}[faithful-diffusion]
% If $\A^m=\{\emptyset , \Omega \}$ for all $m$, faithfulness holds if $\bS$ is a system of linear time-homogeneous diffusion processes.
If the system $\PS^M=(\bA^M,\bX^M)$ is such that $\A^M=\{\emptyset , \Omega \}$ and $\bX^M$ is a linear time-homogeneous diffusion process under $P$, then faithfulness holds for any sequence of nested systems $\bS=(\PS^1, \ldots, \PS^M)$, where $\PS^m=(\bA^m,\bX^m)$ is a system such that $\bX^m\in \bX^M$.
\end{Proposition}

 Even in the true probability the influences for $\bX$ may be
non-causal. However, with the faithfulness assumption two conclusions
can be drawn: (i) if $X_1 \infl {\PS^{m'}} X_2$ then either this
influence is causal or one can find a $\PS^m$, $\PS^{m'} \subset
\PS^m \in \bS$, in which there is a process which influences both
$X_1$ and $X_2$ (a common ancestor in graph terminology): such a
process may be called a confounder in epidemiological terminology;
(ii) if $X_1 \sli {\PS^{m'}} X_2$ this means that  $X_1$ does not
have a causal influence on $X_2$. If an indirect influence in
$\PS^{m'}$ is causal it is stable by considering richer systems
$\PS^m$; direct influences in $\PS^{m'}$ may be related to indirect
causal influences in $\PS^m$.

Now we study criteria of independence of processes, which leads us
to a mathematical proof in our context of the causal interpretation
of a direct influence of a randomized process (our Theorem
\ref{non-influence}).

\begin{Definition}[Dynamical independence]
If  $X_j \sli {\PS^m} X_k$, \mbox{$ X_k \sli {\PS^m} X_j$} and $X_j$
and $X_k$ have no common ancestor, we say that $X_j$ and $X_k$ are
dynamically independent in $\PS^m$.
\end{Definition}

\begin{Definition}[Non-influenced process]
  A process $X_A \in \bX^m$ is a non-influenced process in $\PS^m=(\bA^m,\bX^m)$ for probability $P$ if $X_j \wli {\PS^m} X_A$ for $P$, for all $X_j \in \bX^m$.
\end{Definition}

\begin{Lemma}[Component and group Dynamical
independence]\label{Dynidp-group} If $X_j$ and $X_k$ are dynamically
independent in $\PS^m=(\bA^m,\bX^m)$, it is possible to find $X_A,
X_B, X_C$ such that $\bX^m=(X_A, X_B, X_C)$ where $X_A$ and $X_B$
are non-influenced in $\PS^m$. Conversely any component, say $X_j$,
of $X_A$ and any component of $X_B$, say $X_k$, are dynamically
independent in $\PS^m$.
\end{Lemma}

If $X_C$ is influenced by both $X_A$ and $X_B$, the influence graph
of $\bX^m=(X_A, X_B, X_C)$ where $X_A$ and $X_B$ are non-influenced
 is $X_A \longrightarrow X_C
\longleftarrow X_B$

\begin{Lemma}[Dynamical independence and independence] \label{Dynidp-idp}
 Let $\bX^m=(X_A, X_B, X_C)$. Consider the assumptions : (a) $X_C \wli {\PS^m} X_A$ and $X_C \wli {\PS^m} X_B$,
 (b) $P$ is faithful for the sequence $\bS = (\PS^2,\PS^m)$ with $\PS^2=(\bA,\bX^2)$ with $\bX^2=(X_A,X_B)$ and (c) the decomposition of $(X_A,X_B)$ in $\PS^m$ is that of a diffusion with jumps such that given $\A$ the corresponding SDE satisfies the
 uniqueness conditions in law.

 Consider the two following
propositions :
\begin{itemize}
\item[(i)]$X_A$ and $X_B$ are  independent conditionally on $\A$;
\item[(ii)]$X_B \sli {\PS^m} X_A$, $X_A \sli {\PS^m} X_B$ and $X_{A0}\ind_{\A}X_{B0}$.
\end{itemize} 
Then, under assumptions (a) and (b), (i) implies (ii). Moreover if  (c) holds, the converse is true.
\end{Lemma}

{\bf Remark.} Diffusion with jumps and conditions of uniqueness are given in Jacod
and Shiryaev (§III.2). In proposition (ii) $X_{A0}$ and $X_{B0}$ are the initial values of $X_{A}$ and $X_{B}$.

 Although the Lemma may seem intuitively obvious a general
proof is not simple to find. See Appendix B for an outline of proof.

\begin{Definition}[$\bS$-Non-influenced process]
  Let a system $\PS^{1}$ belonging to a  sequence $\bS$; a process $I \in \bX^1$ is a $\bS$-non-influenced process for probability $P$ if whatever $\PS^m \subset \bS$, $I \ind \A^m$ and  $X^m_j \wli {\PS^m} I$ for $P$, for all $j$.
\end{Definition}

 The only clearly non-influenced processes for $P^*$ are randomised processes, generally randomised attribution of a treatment. In observational studies, the non-influenced quality will always be an assumption. For instance a genetic factor may in some circumstances be considered as non-influenced (Didelez and Sheenan, 2007). However, in our approach genetic factors would generally be considered as (observed or non-observed) attributes, and not part of the state.

 \begin{Theorem}[Non-influence and causality] \label{non-influence}
 Let $\bS$ an approximating sequence for $\PS^{1}=(\bA^1,(I,X_j))$. Suppose that $P^*$ is faithful for any sequence in the associated perfect system
 $\PS^M$, that  $(I, X_j)$ satisfies the assumption (a) and (c) of Lemma \ref{Dynidp-idp} and $I_0
\ind_{\A^m}X_{j0}$ for all $m$. If $I$ is a $\bS$-non-influenced
process for $P^*$ and $I\dinf {\PS^1} X_j$ for $P^*$, then $I$
causally influences $X_j$.
\end{Theorem}
{\bf Proof.} If $I$ did not causally influence $X_j$, we would have
$I \sli {\PS^M} X_j$ for $P^*$. Since $I$ is a non-influenced
process, according to Lemma \ref{Dynidp-group} and Lemma
\ref{Dynidp-idp}, $I\ind_{\mathcal{A}^M}X_j$ and using the fact that
$I \ind \mathcal{A}^m$ for all $m$ it implies
$I\ind_{\mathcal{A}^1}X_j$, and in particular that $I \wli {\PS^1}
X_j$, in contradiction with our assumption. Hence the Theorem.

It is interesting to give a version of the idea of instrumental
variables (Stock, 2001; Angrist et al., 1996; Greenland, 2000)
applied to our context; here the idea is applied only to assess the
causal nature of an influence, while it is often used to estimate
the magnitude of the causal influence in specific models. We have
the following result:

\begin{Lemma}[Instrumental processes] \label{Intrumental}
Under the assumptions of Theorem \ref{non-influence}, if $I$ is a
$\bS$-non-influenced process, $I\dinf {\PS^1} X_k$, and $X_j$ blocks
the paths from $I$ to $X_k$ in system $\PS^M$, then $X_j$ causally
influences $X_k$.
\end{Lemma}
{\bf Proof.} By Theorem \ref{non-influence} we have $I\dinf
{\PS^1} X_k \Longrightarrow I\infl {\PS^M} X_k $. If $X_j$ blocks
the paths from $I$ to $X_k$ in $\bX^M$, then by Lemma
\ref{decomp-infl} we have $X_j\infl {\PS^M} X_k $; hence the
Lemma.

%\begin{Lemma}[No unmeasured confounders] \label{no-unmeasured-confounder}
% Consider $\PS^1=(\bA^1, \bX^1)$ with $\bX^1=(X_j,X_k, X_C)$. If the following assumptions hold under $P*$:  
%(i) if $X_{C'} \dinf {\PS^M} X_j$,  $X_C$ blocks the path from $X_{C'}$ to $X_k$; if $X_{C''} \dinf {\PS^M} X_k$,  $X_C$ blocks the path from $X_{C''}$ to $X_j$;
%(ii)  $X_k \sli {\PS^M} X_j$;
%(iii) $X_{j0} \ind {\A_1} X_{j0}$;
%(iv) $\bA^1=\bA^M$;

% then, if $X_j\dinf {\PS^1} X_k$, $X_j$ causally influences $X_k$.
%\end{Lemma}

%{\bf Proof.} If $X_j$ did not causally influence $X_k$ we would
%have  $X_j \sli {\PS^M} X_k$; with the conditions stated in the
%Lemma, $X_j$ and $X_k$ would then be dynamically independent. With the other assumptions, in
%virtue of Lemma \ref{Dynidp-idp} they would also be dynamically
%independent in $\PS^1$, which contradicts the assumption; hence the
%Lemma.

\subsection{Implications for physics, system biology and epidemiology}

\subsubsection{Physics}

Let us consider as an example the level of the trajectories of planets in the solar system. The physical system
is the set of planets simplified to points in three-dimensional space and we are interested only in their
trajectories. 
%Note that when we consider this system we neglect the influence of events occurring in other
%stellar systems as well as that of events occurring on the planets themselves (such as the flight of
%butterflies on the Earth). 
The state of the system can be represented by a multivariate process $\bX$, the components of
which are the positions and the speeds of the planets in a given set of axes and this process obeys a
differential equation of the type $d\bX_t=g(\bX_t)dt$, where $g(.)$ is a function derived from Newton's law of
mass attraction. $\bX$ is a time-homogeneous Markov process, although degenerated because deterministic. There do
not seem to be processes that can be manipulated in this system. However we believe that the influence of a
planet on the trajectory of another planet may be considered as being of causal nature.

A first instance of the application of physical laws is to predict or to control the state of the system: for instance one can predict eclipses or control the trajectory of a space vessel. In this case we assume that we know the physical law and that we have a good system.

A second instance is that there is a discrepancy between $P^*$ and $P^{\PS}$ for the chosen system. If there is
not much doubt about the physical laws we are applying (here Newton's laws) then it may be deduced that the
system considered is not satisfactory and that it must be increased. A famous example of such an instance is the
discrepancy which appeared between the computed and the observed trajectories of Uranus. Leverrier made
computations which lead to the discovery of Neptune in 1846. He assumed that the discrepancy in the observed
trajectory of Uranus with respect to what was computed using Newton's laws was due to the presence of another
planet: he gave the computed position of this planet to Johann Galle and Louis D'Arrest who found it.

A third instance occurs if in spite of refining the system, a discrepancy persists. Then the physical laws may
be cast into doubt.

% This happened when the 1919 sun eclipse was observed: the conventional theory did not predict the deflection of light rays due to gravitation. This lead to abandon the conventional theory (for these phenomena) in favour of Einstein's theory of relativity.

\subsubsection{Systems biology}
  The model is constructed with partially known mechanisms but some of the influences are unknown and even when causal influences are assumed, their precise forms are unknown. These models can be used to test whether some causal influences exist or to quantify them when they are assumed to exist. We will develop the analysis of the interaction between HIV and the immune system in section \ref{HIV}.

\subsubsection{Epidemiology}

 Most epidemiological studies endeavour to test the influence of a single factor on a disease process. The physical system contains all biological phenomena implied in the disease as well as the factor of interest; in general there is no physical law, only biological plausibility of some causal influences. A typical system is $\bX=(F, D, C)$ where $F$ is the factor of interest, $D$ represents the disease and $C$ are other processes taken from the system. Such a problem is most often modelled with random variables rather than with stochastic processes. The stochastic process framework allows to take into account the dynamics of the phenomena: typically $D$ would be a counting process and the exposure factor $F$ may also vary in time, as is most often the case in reality. The interest often lies in the possible causal influence of $F$ on $D$. Testing whether  $F \dinf {\bX} D$ is generally expressed by saying that we test whether $F$ is a risk factor for $D$ by an analysis adjusted on $C$. 
 It should be possible to formmalize in our framework the condition of  ``no unmeasured confounders''  which makes it possible to conclude that $F$ causally influences $D$. This however requires further work.

 In many simple clinical trials the main interest lies in a particular influence, that of a drug on a clinical endpoint. The aim is to test whether there is a causal influence without trying to understand which basic causal mechanisms may explain it, even if there is a biological plausibility that a certain molecule (or treatment) may have a causal influence on the clinical endpoint considered. That is, most often, we do not have physical laws.  This is why randomised trials have been developed. If $F$ is a treatment that can be randomised, it becomes $\bS$-non-influenced. Then by Theorem \ref{non-influence} it is sufficient to look at the influence of $F$ on $D$ in any model to deduce the presence or absence of causal influence.

\section{Model for the observations}
\label{observations}
In most applications we do not have precise physical laws. Instead of a unique probability we use  a model, that is a family of probability $(P^{\PS^{m}}_{\theta})_{\theta \in \Theta}$ on $\FF^m$. The choice of the model may include scientific knowledge, that is a model can be considered as an incompletely specified physical law. If the system $\PS^{m}$ is rich enough (ideally if it is the ''perfect'' system $\PS^{M}$) and if the knowledge incorporated in the model is correct, the model is well-specified, that is $P^*_{\FF^m} \in (P^{\PS^{m}}_{\theta})_{\theta \in \Theta}$. Even if the model is not well specified it is interesting to find the value $\theta_0$ such that $P^{\PS^{m}}_{\theta_0}$ is the closest to $P^*_{\FF^m}$. Since the latter is unknown we need observations, which by definition are realisations of $\FF^m$-measurable  random variables under probability $P^*_{\FF^m}$. Generally complex systems will be observed with complex observations schemes, leading to incomplete (or coarsened) or indirect observations. Generalising the approach of Heitjan and Rubin (1991) to stochastic processes we may say that the observation, represented by the sigma-field $\Ob^m$, are generated by $g(\bX,G)$, where $G$ is a component which may be deterministic or stochastic. If $G$ is deterministic we have $\Ob^m \subset \FF^m$; if however $G$ is random, $\Ob^m$ is  not a subset of $\FF^m$.

To choose a probability in the model close to $P^*_{\FF^m}$ we must construct an estimator $\hat \theta (\Ob^m)$. For maximum likelihood or maximum penalised likelihood estimators we must compute the likelihood for the observation, which is the Radon-Nykodim derivative of $P^{\PS^{m}}_{\theta}$ relative to a reference probability $P_0$ on the sigma-field $\Ob^m$, and we denote it $\LL_{\Ob^m}^{P^{\PS^{m}}}$. If the mechanism leading to incomplete data (m.l.i.d.) is deterministic this is equal to $\EPZ(\LL_{\FF^m}^{P^{\PS^{m}}}|\Ob^m)$ and this is relatively easy to compute. If not, the issue of ignorability of the m.l.i.d. arises: if the m.l.i.d. is ignorable we can proceed as if it was deterministic and obtain nevertheless the correct inference. For instance Commenges and G\'egout-Petit (2007) computed the likelihood for counting processes observed with a complex, but ignorable, observation scheme. If the m.l.i.d. is not ignorable we have to include $G$ in the system and consider $\bX^{m'}=(\bX^m,G)$; we have then by definition $\Ob^{m'} \subset \FF^{m'}$ and we can apply the above formula. The price to be paid is that we need additional assumptions and the computation of the likelihood may not be easy.

In epidemiology one generally has a sample of observations for a sample of systems indexed by $i$, $i=1,\ldots,n$. The most common framework is that the observations are independently identically distributed. In this framework, if we can describe the system and its observation for a generic item, we can do it for the sample; this is why in this paper we always omit the subscript $i$.

\section{Dynamical models for HIV/AIDS}
\label{HIV}
\subsection{The problematic of AIDS through dynamic influence graphs}
AIDS was identified in 1981 as a life-threatening disease due to acquired immunodeficiency. It was found
that this immunodeficiency was essentially due to a decrease of the number of CD4+ T-lymphocytes. In 1983 it was found
that this decrease was mainly due to the destructive replication of a virus in CD4+ lymphocytes and this virus
was denominated HIV. Thus we can formulate the causal pathway: ``presence of HIV causes low CD4 counts which
causes AIDS which causes death''. Although most researchers would agree with this phrase and think that what is
behind the word ``cause'' are particular biological mechanisms which could be further reduced to biochemical
laws, it remains vague because i) time is only implicitly involved through the fact that cause precedes effect;
ii) each modality is relative to another modality (presence vs absence, low vs high and so on).

The dynamical model approach allows us to make the causal statement more precise. First we construct the processes $I=(I_t$), $T=(T_t$), $A=(A_t$), $D=(D_t$): $I$ is a counting process representing HIV infection, $T$  has a continuous state-space and represents CD4+ T-lymphocytes count; $A$ and $D$ are counting processes representing AIDS and death respectively. We can express the causal structure by the influence graph:
$$ I \dinf{} T \dinf{} A \dinf{} D$$
Indeed we know from the results of research (involving virology, immunology, and clinical research) that these influences can be interpreted as causal. It is interesting to note that we consider that $ I \dinf{} T$ is causal although it is difficult to manipulate $I$. A more detailed description of the infection can be made by introducing the viral load process $V=(V_t)$. There is of course a direct influence of $I$ on $V$ because if $I_t=0$ then $dV_t=0$. When considering the evolution of infected subjects the process of interest is $V$ (not $I$ which is identically equal to one in these patients).
\subsection{From descriptive to mechanistic models}
In the conventional epidemiological and biostatistical literature, linear mixed-effect models have been used to
analyse separately repeated measurements of CD4 counts and viral load. For instance to analyse viral load
following initiation of a highly active anti-retroviral therapy (HAART) a linear-mixed effect model with two
slopes has been used (Jacqmin-Gadda et al., 2000). Potential observations $Y_{j}$ are the viral load at time
$t_{j}$, or a logarithmic transformation of the viral load; for simplicity we will ignore these normalising
transformations here. Some data may be missing (a non-ignorable mechanism here): $Y_{j}$ was observed only if
$Y_{j} > \eta$, where $\eta$ is a detection limit, while $1_{Y_{j}>\eta}$ was always observed. The model can be
written as:
\begin{equation} \label{mixlin1} Y_{j}=\beta_0+ a_{0} +(\beta_1+ a_{1} +\gamma_{1}A) \min (t_{j},t_*) +(\beta_2+a_{2}+\gamma_{2}A)(t_{j}-t_*) I_{t_{j}>t_*}+ \varepsilon_{j},\end{equation}
where $\beta_{0}$, $\beta_{1}$, $\beta_{2}$, are parameters for the intercept, first and second slopes respectively and $a_{0}$, $a_{1}$, $a_{2}$, are independent normal random effects on the intercept, first and second slopes respectively; $t_*$ is the time of change of slope (supposed known), $A$ indicates the treatment and $\varepsilon_{j}$ are normal variables with zero expectations: they may be independent or have a correlation structure.
In the dynamical model representation, this model can be written in terms of the process $V=(V_t)$ living in continuous time, representing the concentration of virus at time $t$. There are at least two ways of representing the random effects: they could be degenerate components of the state or they could be random attributes. We adopt the latter which leads to the simplest expression:
\begin{equation} \label{eqdif1} dV_t=[(\beta'_1+\gamma_1 A) I_{t\le t_*}+(\beta'_2+\gamma_2 A)I_{t > t_*}]dt +\sigma dW_t,\mbox{ with } Z_0=\beta'_0\end{equation}
where $\beta'_0$ is a random initial condition and $\beta'_1$ and $\beta'_2$ are considered as random
attributes; the link with the above model is that $\beta'_j$ has expectation $\beta_j$ and variance ${\rm var}~
a_j$. The observation (treating the observation times as fixed) is $\Ob=\sigma
(1_{Y_{j>\eta}},1_{Y_{j>\eta}}Y_j, j=1,\ldots,m)$ where $Y_{j}=V_{t_{j}}+\varepsilon'_{j}$. Note that the error
$\varepsilon_{j}$ of model (\ref{mixlin1}) is the sum of the value of the martingale at $t_j$ and the observation
error in  model (\ref{eqdif1}): we have $\varepsilon_{j}=W_{t_j}+\varepsilon'_{j}$. The models for the observations may be
the same if the correlation structure of the $\varepsilon_{j}$ in model (\ref{mixlin1}) is compatible with that
produced by model (\ref{eqdif1}). The graph of this process is not very interesting since only $A$ influences $Z$.

A more elaborate model was proposed by Thi\'ebaut et al. (2005). This was a multivariate linear mixed model for jointly modelling viral load and CD4, together with a possibly informative drop-out. For each of the two markers there were two slopes with a fixed and a random effect (as in the previous model). We leave aside here a certain number of features of that paper, including modelling of the drop-out and of explanatory variables, to focus on how the link between observations of HIV load and CD4 counts was modelled. The model can be written:
$$Y^1_{j}=\beta^1_0+ a^1_{0} +(\beta^1_1+ a^1_{1} ) \min (t_{j},t_*) +(\beta^1_2+a^1_{2})(t_{j}-t_*) I_{t_{j}>t_*}+ \varepsilon^1_{j},$$
$$Y^2_{j}=\beta^2_0+ a^2_{0} +(\beta^2_1+ a^2_{1} ) \min (t_{j},t_*) +(\beta^2_2+a^2_{2})(t_{j}-t_*) I_{t_{j}>t_*}+ \varepsilon^2_{j}.$$ where $\varepsilon^1_{j}$ and $\varepsilon^2_{j}$ are zero expectation normal variables. For fixed $j$, $\varepsilon^1_{j}$ and $\varepsilon^2_{j}$ are independent; the sequences $\varepsilon^k_{j}, j=1,\ldots, m$ for $k=1,2$ may be formed of independent variables or have a correlation structure. The link between HIV load and CD4 counts was expressed by correlations of the random effects $a^1_l$ and $a^2_l$, $l=0,1,2$. In particular we could expect negative correlations between the slopes of HIV load and CD4 counts, which was indeed observed when fitting the model to the data of a therapeutic trial (better viral response was correlated to better immune response).

The model can be expressed in the dynamical framework as:
$$dV_t=[\beta'^1_1 I_{t\le t_*}+\beta'^1_2I_{t > t_*}]dt+\sigma_1 dW_{1t},\mbox{ with } V_0=\beta'^1_0$$
$$d\bar T_t=[\beta'^2_1 I_{t\le t_*}+\beta'^2_2I_{t > t_*}]dt+\sigma_2 dW_{2t},\mbox{ with } \bar T_0=\beta'^2_0$$
where $V_t$ is the logarithm of the viral load and $\bar T_t$ the CD4 counts at time $t$, $\beta'^k_l=\beta^k_l+ a^k_l, k=1,2; l=1,2$. As in the previous model there are several ways of treating the random effects; for instance we may consider them as random attributes.
The observation is
\begin{equation} \label{Observation} \Ob=\sigma (1_{Y^1_{j>\eta}},1_{Y^1_{j>\eta}}Y^1_j, Y^2_j, j=1,\ldots,m) \mbox { with } Y^1_{j}=V_{t_{j}}+\varepsilon'^1_{j} \mbox { ;  } Y^2_{j}=\bar T_{t_{j}}+\varepsilon'^2_{j} \end{equation}

It is clear from the differential equations above that there is no influence of $V$ on $\bar T$ whatever the values of the parameter: the influence
graph is made of two disconnected vertices. We might have treated the random effects as ancestors, but in this
representation too, there is no direct nor indirect influence of $V$ on $\bar T$. In this model $\bar T$ is SCLI
from $V$ which does not fit with the known mechanism of the infection. So although this model succeeded in fitting
the data better than separate linear mixed models, it is unable to capture any relevant causal influence.

There are different models in which we can express that viral load influences CD4. Having made a clear distinction between the ``model for the system'' and the ``model for the observation'' it is natural to construct a model including components that are not observed at all, but which will be more satisfying with respect to the way it represents the biological mechanisms. One may distinguish infected and un-infected cells and take into account the causal influences in the ODE system (Ho et al., 1995; Perelson et al., 1996). Still a more satisfying model distinguishes between quiescent ($Q$)and activated CD4 ($T$) and between infectious ($V_I$) and non-infectious ($V_{NI}$) virus. 
Note that distinguishing quiescent and activated CD4 is a way of enriching the state without simply adding a new component.
To write the differential equation for the model one uses additional assumptions which are plausible in view of the knowledge of the biological mechanisms: for instance we assume that new CD4+ T lymphocytes are produced (by the thymus) at a rate $\lambda$, that only activated cells can be infected, that the probability of meeting of a cell and a virion is proportional to the product of their concentrations. The model proposed by  Guedj, Thi\'ebaut and Commenges (2007) was:

\begin{eqnarray*}\label{(1)}
dQ_t&=& (\lambda + \rho T_t - \alpha Q_t- \mu_{Q} Q_t)dt \\
dT_t&=& (\alpha Q_t  - (1-\eta \1{I^{RT}_t=1}  )\gamma T_tV_{It} - \rho T_t - \mu_{T}T_t)dt \\
dT^{*}_t&= &[(1-\eta \1{I^{RT}_t=1}  )\gamma  T_tV_{It} - \mu _{T^{*}} T^{*}_t]dt \\
dV_{It}& =& (\omega \mu _{T^{*}_t}\pi T^{*}_t - \mu _{v}V_{It})dt \\
dV_{NIt}& =& [(1-\omega)\mu _{T^{*}_t} \pi T^{*}_t - \mu _{v}V_{NIt}]dt \\
\end{eqnarray*}
where $I^{RT}$ is the process indicating whether a treatment based on an inhibitor of the reverse transcriptase
is taken at time $t$. If we consider the framework of a controlled clinical trial this process is non-influenced
and controlled (because its trajectory is obtained by randomisation).
% In principle we should distinguish between
%attribution of treatment (which is really controlled) uptake of the drug (which is modulated by compliance) and
%concentration of the drug (which obeys pharmacokinetics laws): for the sake of simplicity we do not develop these
%distinctions here. 
Guedj, Thi\'ebaut and Commenges (2007) assumed that some
parameters were random. Such parameters may be considered as random attributes while fixed parameters may be considered as constants of a
``physical law''. Note that the system is time-homogeneous, which is satisfactory from an explanatory point of
view. Moreover, as we noted in section \ref{detcase}, this makes it possibile to draw the influence graph of a
deterministic model. We could also consider a stochastic differential equation system but inference in this
context is very challenging. The observation  is the same as in (\ref{Observation}), with
$V_{t_j}=V_{It_j}+V_{NIt_j}$ and $\bar T_{t_j}= Q_{t_j}+T_{t_j}+T^*_{t_j}$.

\begin{figure}[h]

\centering

\includegraphics[scale=0.6]{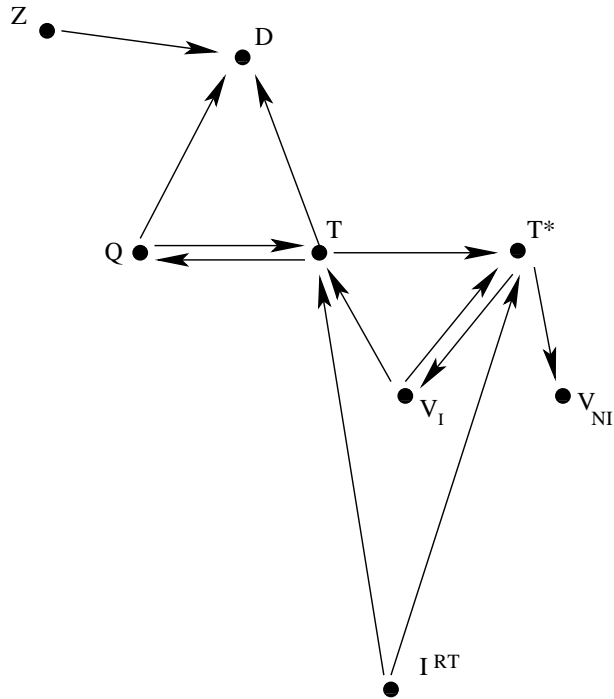}

\caption{Graph for the mechanistic HIV model.\label{HIV-complex-RTI}}
\end{figure}

We could consider mixing this model for the markers with a model for an event such as an opportunistic disease,
adding the component $D=(D_t)$ which is a counting process. The risk of the opportunistic disease may be
considered as depending on the concentration of CD4+ T lymphocytes, so that keeping the framework of a
time-homogeneous Markov model we can propose a proportional hazard model (but with constant base-line risk
$\gamma$):
$$dD_t=I_{\{D_{t-}=0\}}\gamma \exp (\beta_1 Q_t + \beta_2 T_t + \beta_3 Z)+dM_t,$$

where $Z$ is an explanatory variable.
The graph for such a model is given in Figure \ref{HIV-complex-RTI}.

\begin{figure}[h]

\centering

\includegraphics[scale=0.6]{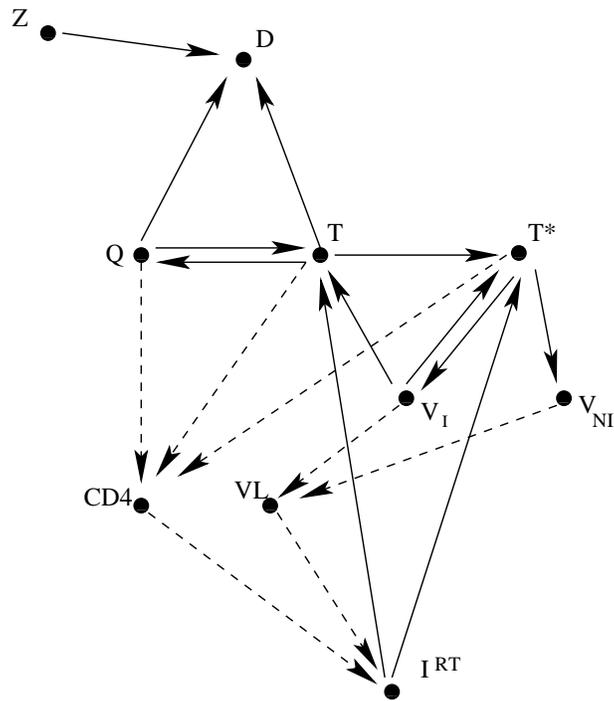}

\caption{Graph for the mechanistic HIV model.\label{HIV-complex-RTI-doctor}}
\end{figure}

Note that if the treatment was an inhibitor of protease the graph would be different: the inhibitor of protease
influences $V_I$ and $V_{NI}$. Also, in an observational study, the treatment is in fact influenced by the
information on the clinical and biological state of the patient. If we want to represent this situation we
have to include the medical doctor in the system: the doctor may decide to modify the treatment after having
been informed of the measurement of  viral load (VL) and of CD4 counts (CD4); note that the processes VL and CD4
and are different from $V$ and $\bar T$ because they carry the information on measurements of these
processes, that is $(V_t, \bar T_t)$ carry the observation contained in $\Ob$ (see (\ref{Observation})) up to
time $t$.  Then the graph could be as shown in Figure \ref{HIV-complex-RTI-doctor}, where we have represented
by dotted lines the influences of the marker processes on their measurements and the influence of these
measurements on the treatment, through the decision of the doctor.

\section{References}

\setlength{\parindent}{0.0in}

\noindent Aalen, O. (1987) Dynamic modelling and causality. {\em Scand. Actuarial J.}, {\bf 1987}, 177-190.\vspace{3mm}

 Aalen, O. and Frigessi, A. (2007) What can statistics contribute
   to a causal understanding? {\em Scand. J. Statist.}, {\bf 34}, 155-168.\vspace{3mm}

%\noindent Andersen, P.K., Borgan \O., Gill R.D.  \&  Keiding N. (1993) {\it Statistical Models Based on Counting Processes}. New-York: Springer-Verlag.\vspace{3mm}

Angrist, J., Imbens, G., and Rubin, D. (1996) Identification of causal effects using
instrumental variables. {\em J. Am. Statist. Ass.}, {\bf 91}, 444-455.\vspace{3mm}

Batterman, R.W. (2002) {\em The devil in the details}. Oxford University Press. \vspace{3mm}

Bunge, M. (1979) {\em Causality and modern science}. Dover Publications, New-York.\vspace{3mm}

\noindent Commenges, D. and G\'egout-Petit, A. (2007) Likelihood for generally coarsened observations from multi-state or counting process models. {\em Scand. J. Statist.}, {\bf 34}, 432-450.\vspace{3mm}

Commenges, D., Joly, P., G\'egout-Petit, A. and Liquet, B. (2007). Choice between semi-parametric estimators of Markov and non-Markov multi-state models from generally coarsened observations. {\em Scand. J. of Statist.}, {\bf 34}, 33-52.\vspace{3mm}

Cox, D. and Wermuth, N. (1996) {\em Multivariate Dependencies: Models, Analysis and Interpretation}. Chapman \& Hall, London.\vspace{3mm}

Dawid, A. P. (1979) Conditional independence in statistical theory (with Discussion).
{\em J. Roy. Statist. Soc. A}, {\bf 41}, 1-31.\vspace{3mm}

Dawid, A. P. (2000) Causal inference without counterfactuals. {\em  J. Am. Statist. Ass.}, {\bf 95}, 407-448.\vspace{3mm}

%Didelez, V. (2000) {\em Graphical Models for Event History Data Based on Local Independence}, PhD thesis, Universitat Dortmund, Germany.\vspace{3mm}

%Didelez, V. (2005) Mendelian randomisation and instrumental variables: What can and what can't be done (with N. Sheehan), Research Report 05-02, Department of Health Sciences, University of Leicester.\vspace{3mm}

Didelez, V. (2007). Graphical models for composable finite Markov processes, {\em Scand. J. Statist.}, {\bf 34}, 169-185. \vspace{3mm}

Didelez, V. and Sheehan,  N.A.(2007) Mendelian randomisation as an instrumental variable approach to causal inference. {\em Statistical Methods in Medical Research}, in press.

%Didelez, V. (2007) Mendelian randomisation: why epidemiology needs a formal language for causality, to appear in: F. Russo and J. Williamson (eds.), Causality and Probability in the Sciences, College Publications London, 2007.\vspace{3mm}

Eerola, M. (1994). {\em Probabilistic Causality in Longitudinal Studies}. Lecture Notes in Statistics, Springer, New-York.\vspace{3mm}

Fosen, J., Ferkingstad, E., Borgan, O and Aalen, O.O. (2006) Dynamic path analysis: a new approach to analyzing time-dependent covariates. {\em Lifetime Data Analysis}, {\bf 12}, 143-167.  \vspace{3mm}

%Ganiayre, J., Commenges, D. and Letenneur, L. (2007) A Latent Process Model for Dementia and Psychometric Tests. {\em Lifetime Data Analysis}, in revision.\vspace{3mm}

Geneletti, S. (2007) Identifying direct and indirect effects in a non-counterfactual framework. {\em J. Roy. Statist. Soc. A}, {\bf 69}, 199-215.\vspace{3mm}

Gibbs, A.L. and Su F.E. On Choosing and Bounding Probability Metrics. {\em Int Stat. Rev.} {\bf 70}, 419-435.\vspace{3mm}

Gill, R. D., van der Laan, M. J. and Robins, J.M. (1997). Coarsening at random: characterizations, conjectures
and counter-examples, in: {\em State of the Art in Survival Analysis}, D.-Y. Lin and T.R. Fleming
(eds), Springer Lecture Notes in Statistics 123, 255-294.  \vspace{3mm}

Granger, C. W. J. (1969) Investigating causal relations by econometric
models and cross-spectral methods. {\em Econometrica} {\bf 37}, 424-438.\vspace{3mm}

Greenland, S. (2000) An introduction to instrumental variables for epidemiologists. {\em Int. J. Epidemiol.}, {\bf 29}, 722-729.\vspace{3mm}

Guedj, J., Thi\'ebaut, R. and Commenges, D. (2007) Maximum Likelihood Estimation in Dynamical Models of HIV. {\em Biometrics}, in press.\vspace{3mm}

Ho, D.D., Neumann, A.U., Perelson, A.S., Chen, W., Leonard, J.M. and Markowitz, M.(1995) Rapid turnover of plasma virions and CD4 lymphocytes in HIV-1 infection. {\em Nature} 1995, 373(6510), 123-126.\vspace{3mm}

Holland, P. (1986) Statistics and Causal Influence. {\em J. Am. Statist. Ass.}, {\bf 81}, 945-960.\vspace{3mm}

Jacqmin-Gadda, H., Thi\'ebaut, R., Ch\^ene, G. and Commenges, D. (2000) Analysis of left-censored longitudinal data with application to viral load in HIV infection. {\em Biostatistics}, {\bf 1}, 355-368.\vspace{3mm}

%\noindent Jacod, J. (1975) Multivariate point processes: predictable projection; Radon-Nikodym derivative, representation of martingales. {\it Z. Wahrsheinlichkeitsth.}, {\bf 31}, 235-253.\vspace{3mm}

Jacod, J. and Shiryaev, A.N. (1987). {\it Limit Theorems fo Stochastic Processes}. Berlin:
Springer-Verlag.\vspace{3mm}

Jazwinsky, A.H. (1970) {\em Stochastic processes and filtering theory}, Academic Press, New-York.\vspace{3mm}

Joreskog, K. (1978) Structural analysis of covariance and correlation matrices. {\em Psychometrika} {\bf 43}, 443-477.\vspace{3mm}

%Kallenberg, O. (2001) {\em Foundations of modern probabilities}. Springer Verlag, New-York.\vspace{3mm}

%Kalman, R.E. and Bucy, R.S. (1961) New results in linear filtering and prediction theory. {\em J. Basic Eng.}, {\bf 83}, 95-108.\vspace{3mm}

%Lipster, R. S. and Shiryaev, A.N. (2001). {\it Statistics of Random Processes}. New York: Springer-Verlag.\vspace{3mm}

Murrray, C.D.  and Dermott, S.F. (1999) {\em Solar System Dynamics}. Cambridge University Press.\vspace{3mm}

Oksendal, B. (2000) {\em Stochastic Differential Equations}. Springer Verlag, New-York.\vspace{3mm}

Pardoux, E., (1991) Filtrage non lin\'eaire et \'equations aux d\'eriv\'ees partielles
   stochastiques associ\'ees. {\em \'Ecole d'\'Et\'e de Probabilit\'es de Saint-Flour XIX---1989}, Lecture Notes in Math., {\bf 1464}, 67-163,  Springer, Berlin. \vspace{3mm}

Pearl, J. (2000) {\em Causality : Models, Reasoning, and Inference}, Cambridge
University Press.\vspace{3mm}

%Pearl, J. (2001) Direct and indirect effects, in J. Breese and D. Koller, eds,
%Proceedings of the Seventeenth Conference on Uncertainty in Articial
%Intelligence, Morgan Kaufmann Publishers, San Francisco, 411-420.\vspace{3mm}

Perelson, A.S., Neuman, A.U., Markowitch, M., Leonard, J.M. and Ho, D.D. (1996) HIV-1 dynamics in vivo: virion clearance rate, infected cell life-span, and viral generation time. {\em Science}, {\bf 271}, 1582-1586.\vspace{3mm}

Revuz, D. and Yor, M. (1991) {\em  Continuous Martingales and Brownian Motion},  New-York: Springer-Verlag.\vspace{3mm}

Robins, J., Scheines, R., Spirtes, P. and Wasserman, L. (2003) Uniform consistency in causal inference. {\em Biometrika}, {\bf 90}, 491, 515\vspace{3mm}

Rubin, D. B., (1974) Estimating Causal Effects of Treatments in Random-
ized and Nonrandomized Studies. {\em Journal of Educational Psychology}, {\bf 66}, 688-
701.\vspace{3mm}

Salmon, W. (1984) {\em Scientific Explanation and the Causal Structure of the World}, Princeton University Press.\vspace{3mm}

Spirtes, P., Glymour, C. and Scheines, R. (2000) {\em  Causation, Prediction, and Search}, MIT Press, Cambridge, Massachussetts.\vspace{3mm}

Strevens, M. (2005) How Are the Sciences of Complex Systems Possible? {\em Philosophy of Science}, {\bf 72}, 531-556.\vspace{3mm}

Schweder, T. (1970) Composable Markov processes. {\em J. Appl. Probab.}, {\bf 7},
400-410.\vspace{3mm}

Thi\'ebaut, R., Jacqmin-Gadda, H., Babiker, A. and Commenges, D. (2005) Joint Modelling of bivariate
longitudinal data with informative drop-out and left-censoring, with application to the evolution of CD4+ cell
count and HIV RNA viral load in response to treatment of HIV infection. {\em Statist. Med.}, {\bf 24},
65-82.\vspace{3mm}

van der Laan, M. and Robins, J. (2002) {\em Unified methods for censored longitudinal data and causality}, Springer, New-York. \vspace{3mm}

von Bertallanffy, L. (1968). {\em General System Theory}. New-York: George Braziller.\vspace{3mm}

Wright, S. (1921) Correlation and causation. {\em Journal of Agricultural Re-
search}, {\bf 20}, 557-585.\vspace{3mm}

Wright, S. (1934) The method of path coefficients. {\em Ann. Math. Statist.}, {\bf 5}, 161-215.\vspace{3mm}

%Wu, H.L., Ding, A.A., DeGruttola, V. (1998) Estimation of HIV dynamic parameters. {\em Statist. Med.}, {\bf 17}, 2463-2485.\vspace{3mm}

\section*{Appendix A: Faithfulness in diffusion processes}\label{kalman}
We study the faithfulness property in the case of a system of
linear diffusions. For sake of simplicity we consider a process $\bX^3$, $\bX^3_t=(X_{1t}, X_{2t},
X_{3t})$ where $X_1$, $X_2$, $X_3$ are univariate processes. Let us define the processes $\bX^3$ by the following linear stochastic
differential equations with constant coefficients:
\begin{equation}\label{EP}\left\{\begin{array}{ccl}
           dX_{1t} & = & (a_1X_{1t} + b_1X_{2t}+ c_1X_{3t})dt +
            dW_{1t}\\ dX_{2t} & = & (a_2X_{1t} + b_2X_{2t}+ c_2X_{3t})dt +
            dW_{2t}\\ dX_{3t} & = & (a_3X_{1t} + b_3X_{2t}+ c_3X_{3t})dt +
            dW_{3t}\\
         \end{array}\right.
\end{equation}with initial conditions $X_{10}=X_{20}=X_{30}=0$ and $(W_1,W_2,W_3)$ are independent Brownian motions.
We are interested in the semi-martingale decomposition of $\bX^2_t=(X_{1t},
X_{2t})$ in its own filtration $(\X_{1t} \vee \X_{2t})$. If we note
$E[X_{3t}|\X_{1t} \vee \X_{2t}]=\hat{X}_{3t}$ and using the innovation
theorem, we find that:
\begin{equation}\label{TI}\left\{\begin{array}{ccl}
           dX_{1t} & =  &    (a_1X_{1t} + b_1X_{2t}+ c_1\hat{X}_{3t})dt +   \underbrace{dW_{1t} + c_1(X_{3t}-\hat{X}_{3t})dt}_{=dM_{1t}}\\\
            dX_{2t} & =  &    (a_2X_{1t} + b_2X_{2t}+ c_2\hat{X}_{3t})dt +   \underbrace{dW_{2t} + c_2(X_{3t}-\hat{X}_{3t})dt}_{=dM_{2t}}\\
         \end{array}\right.
\end{equation}$M_{1t}$ and $M_{2t}$ are independent Brownian motions in the filtration
$(\X_{1t} \vee  \X_{2t})$.

If we suppose that the coefficient $b_1 \neq 0$ , the probability
would not be faithful, if $X_2\wli {\bX^2} X_1$ that is $b_1X_{2t}+
c_1\hat{X}_{3t}= f(\X_{1t})$. We use the linear filtering equations
given in (Pardoux, 1991):

\begin{eqnarray}
d\hat{X}_{3t} & = & \left[X_{1t} (a_3-R_t(a_1c_1+a_2c_2)) +
X_{2t}(b_3-R_t(b_1c_1+b_2c_2))
 + \hat{X}_{3t}(c_3-R_t(c_1^2+c_2^2))\right]dt \nonumber\\
 &   & +  R_tc_1dX_{1t} +
 R_tc_2dX_{2t} \\
 dR_t & = & (2c_3 R_t +1 - R_t^2(c_1^2+c_2^2))dt \label{Ricat}
\end{eqnarray}A necessary condition in order to delete the dependence of $X^1$
towards $X^2$ is that the part directed by $dW_2$ in $b_1X_{2t}+
c_1\hat{X}_{3t}$ equals 0, that is $( b_1=-R_tc_1c_2)$. If we remark
that $R_t$ which is a solution of the Riccati differential equation
(\ref{Ricat}), cannot be constant, we conclude that the model is
faithful.

Now if we suppose that the coefficients are no longer constant and
are deterministic time functions and if we suppose that the
following relation is true
\begin{equation}\label{hyp1} b_1(t)=-R_tc_1(t)c_2(t) \end{equation}
The part driven by $dX_{2t}$ in $(_1(t)X_{2t}+ c_1(t)\hat{X}_{3t})$
disappears. According to (\ref{hyp1}) and noting $Z_t= b_1(t)X_{2t}+
c_1(t)\hat{X}_{3t}$, (for convenience, we sometimes omit the
dependence of the coefficients
$(b_1(t),b_2(t),b_3(t),c_1(t),c_2(t),c_3(t))$ on $t$)
\begin{eqnarray*}
 dZ_t & = & b_1(t) dX_{2t}+
c_1(t)d\hat{X}_{3t} + (b'_1(t)X_{2t}+ c'_1(t)\hat{X}_{3t})dt \\
& = & c_1(t) \left[X_{1t} (a_3-R_t(a_1c_1+a_2c_2)) +
X_{2t}(b_3-R_t(b_1c_1+b_2c_2))
 + \hat{X}_{3t}(c_3-R_t(c_1^2+c_2^2))\right]dt \\
 &   & \qquad  +  R_tc_1dX_{1t} +
 \left[(b'_1(t)X_{2t}+ c'_1(t)\hat{X}_{3t})\right] dt \\
 & = &  X_{1t} c_1(t)(a_3-R_t(a_1c_1+a_2c_2)) +
X_{2t}\left[c_1(t)(b_3-R_t(b_1c_1+b_2c_2))+b'_1(t)\right]
 \\
 &   & \qquad + \hat{X}_{3t}\left[c_1(t)(c_3-R_t(c_1^2+c_2^2))+ c'_1(t)\right]dt +  R_tc_1dX_{1t}
\end{eqnarray*}
$Z_t$ is the solution of a stochastic differential equation only driven
by $X_{1t}$ if :
 $$c_1(t)\left[c_1(t)(b_3-R_t(b_1c_1+b_2c_2))+b'_1(t)\right] =b_1(t)\left[c_1(t)(c_3-R_t(c_1^2+c_2^2))+
 c'_1(t)\right]$$Using (\ref{hyp1}) to
substitute $b_1(t)$ we can show that if
 $b_3= R_t(b_2c_2+c'_2+c_2-c_2c_3) +R_{2t}c_2^3$, $Z_t$ is driven by
 $X_{1t} $ and
 the property of faithfulness falls.

This case is extreme. In fact if it holds, the dynamic of $b_1(t)$,
$b_2(t)$ and $b_3(t)$ is imposed by those of $c_1(t)$, $c_2(t)$
and $c_3(t)$

\section*{Appendix B: Proof of Lemma \ref{Dynidp-idp}}\label{prooflemmaidp}
{\bf Proof:}

  Let us first prove that (i) implies (ii). Consider the
Doob-Meyer decomposition of $X_A$ in the filtration ${\cal A} \vee
\X_{At}$: $X_{At}= \Lambda_{At} + M_{At}$. By (i), we have $E[M_{At}
-M_{As}|{\cal A} \vee \X_{As} \vee \X_{Bs}] = E[M_{At} -M_{As}|{\cal A}
\vee \X_{As} ]$ and thus the Doob-Meyer decomposition of $X_A$ is the
same in the filtrations $({\cal A} \vee \X_{At})$ and $({\cal A} \vee \X_{At} \vee
\X_{Bt})$. This implies $X_B \sli{\PS^2} X_A$ .
By symmetry, we have $X_A \sli{\PS^2} X_B$ and
(ii) follows in $\PS^2$. Now by the faithfulness property, we have (ii)
in all system $\PS^m$ with $\PS^2 \subset \PS^m$.

As for the converse, we prove it in the case of a process satisfying
a SDE governed by a Brownian motion in $(\FF_t)$ with  $\FF_t={\cal A} \vee \X_{At}
\vee  \X_{Bt}$:
\begin{eqnarray}
    X_{At} & = &  X_{A0} + \int_0^tf(\X_{As}, \alpha_0)ds + \int_0^t\sigma_{As}dW_{As} \label{xa}\\
    X_{Bt} & = & X_{B0} +\int_0^tg(\X_{Bs},\beta_0)ds + \int_0^t\sigma_{Bs}dW_{Bs} \label{xb}
  \end{eqnarray}where $(\alpha_0, \beta_0)$ is ${\cal A}$-measurable,
 $\sigma_A$ and $\sigma_B$ are deterministic ({\bf A2}) and
$W_A$ and $W_B$ are two independent Brownian motions ({\bf A1}). We
suppose that given $(\alpha_0, \beta_0)\in {\cal A}$,
   the SDE satisfies assumption assuring uniqueness in law (see for
instance Revuz and Yor, 1991: Definitions IX.1.3 and IX.1.4
and Corollary IX.1.14 for the conditions). As by assumption, $X_A$
and $X_B$ are non-influenced in $({\cal A},\mathbf{X}^m)$, then
whatever the system  $\PS^{m'}=(\bA,\bX^{m'})$ such as $\bX^2 \subseteq
\mathbf{X}^{m'}\subseteq \mathbf{X}^{m}$ the process $(X^A,X^B)$
always satisfies  the same SDE.

However, we can take  a new probability space $(\Omega',\FF')$
endowed with two independent Brownian motions $W_{A'}$ and $W_{B'}$
 and construct two independent processes $X_{A't}-X_{A'0}$ and $X_{B't}-X_{B'0}$
on it with $(X_{A'0},X_{B'0}) =^{\cal{L}} (X_{A0},X_{B0})$ and with
$X^{A'}$ satisfying SDE (\ref{xa}) driven by  $ W_{A't} $ and
$X^{B'}$ satisfying SDE (\ref{xb}) driven by  $ W_{B't} $. By the
first part of this demonstration, the decomposition of
$(X_{A'},X_{B'})$ in $\X_{A'} \vee \X_{B'}$ is given by the joint
system of the 2 equations satisfied by $X_{A'}$ and $X_{B'}$ in her
own filtration. The vector $(X_A,X_B)$ and $(X_{A'},X_{B'})$
satisfies the same SDE, by uniqueness in law this implies the
conditional independence between $X_A$ and $X_B$ given ${\cal A}$.

 Using the same reasoning one can extend the result to any
 diffusion system with jumps (defined in Jacod and
Shiryaev, 1987: p. 155) satisfying the condition of uniqueness in law given
${\cal  A}$ (see theorem III.2.32 in Jacod and Shiryaev, 1987).

\end{document}